\documentclass[12pt,a4paper]{article}
\usepackage[english]{babel}
\usepackage[utf8x]{inputenc}
\usepackage{amsmath}
\usepackage[margin=1in]{geometry}
\usepackage{graphicx,cite}
\usepackage[section]{placeins}
\usepackage{float}
\usepackage{ amssymb }
\usepackage[svgnames]{xcolor}
\usepackage[colorlinks=true, linkcolor=Maroon, urlcolor=Black]{hyperref}
\usepackage{amsthm}
\usepackage{amsfonts}
\usepackage{epsfig}
\usepackage{ textcomp }
\usepackage{mathtools}
\usepackage{authblk}
\usepackage{etoolbox}
\usepackage[bottom]{footmisc}
\usepackage{lipsum}

\newcommand{\bb}{\begin{equation}}
\newcommand{\ee}{\end{equation}}

 \newtheorem{thm}{Theorem}
 
 \newtheorem{prop}[thm]{Proposition}
 \newtheorem{lem}[thm]{Lemma}

\title{3-Cube-Free Sets in Cyclic Groups}
\date{(November 2026)}
\author{\v{Z}arko Ran\dj elovi\'c*}

\begin{document}
\maketitle

\begin{abstract}
A set is called \emph{$3$-cube-free} if it contains no configuration of the following form: $\{x,y,z,x+y,x+z,y+z,x+y+z\}$. We study the maximum density of $3$-cube-free sets in finite cyclic groups. As observed by Meng, every 3-cube-free subset of \(\mathbb Z_m\) has size at most \(2m/3\), and this bound is sharp whenever \(3\mid m\). Our main result shows that there exists an absolute constant $\varepsilon>0$ such that, for every $m$ not divisible by $3$, every $3$-cube-free subset of $\mathbb{Z}_m$ has size at most $\left(2/3-\varepsilon\right)m$. Thus, in this setting, the natural $2/3$ barrier can be uniformly improved. Long and Wagner conjectured that every 3-cube-free subset of \(\mathbb Z_{2^n}\) has size at most \(5\cdot 2^n/8\). We also establish this conjecture for 3-cube-free sets whose density in the layer of odd residues is at least \(2/3\). Together with induction, this reduces the conjecture to the case where this density lies in $(5/8,2/3)$.
\end{abstract}

\section{Introduction}

For a positive integer $d$, a \emph{projective $d$-cube} is a set of
the form
\[
 B_d(x_1,\ldots,x_d)
 =
 \left\{
 \sum_{i\in I}x_i \mid
 \emptyset\neq I\subseteq[d]
 \right\}.
\]
A set containing no projective $d$-cube is said to be
\emph{$d$-cube-free}. In particular, a projective $3$-cube is a
configuration 
\[
 \{x,y,z,x+y,x+z,y+z,x+y+z\}.
\]
Questions concerning monochromatic cubes and cube-free sets have a
long history in additive and extremal combinatorics. Rado
\cite{Rado1933} showed that
for every fixed number of colors $c$ and every fixed $d$, there is a 
sufficiently large $n$ such that whenever $[n]$ is colored in $c$ colors there is a monochromatic projective
$d$-cube. Projective
cubes themselves form a natural higher-dimensional generalization of Schur
triples, the case $d=2$. For results on sum-free sets in finite abelian groups,
see \cite{GreenRuzsa2005}.

More recently, Long and Wagner
\cite{LongWagner2019} studied projective cube-free subsets of
$\mathbb Z_{2^n}$. Among other results, they determined the largest
$d$-cube-free subsets when $d$ is a power of two, revealing a close
connection between the extremal structure of cube-free sets and the
natural $2$-adic layers of $\mathbb Z_{2^n}$.

The case of $3$-cubes exhibits rather different behaviour. A recent
paper of Meng \cite{Meng2026} studies $d$-cube-free subsets of
$\mathbb Z_N$ and proves, in particular, that every $3$-cube-free
subset of $\mathbb Z_N$ has size at most $2N/3$. When $3\mid N$, this
is sharp, as witnessed by the set
\[
 \{x\in\mathbb Z_N\mid x\not\equiv0\pmod 3\}.
\]
The constant $2/3$ also arises naturally in the corresponding problem
for subsets of $[N]$. Eberhard and Pohoata conjectured (see
\cite{Meng2026,EberhardMellon2025}) that every $3$-cube-free subset of $[N]$ has size at
most $(2/3+o(1))N$. There are several different constructions
showing that the constant $2/3$ would be best possible. Besides the
sets of integers not divisible by $3$ and the interval $(N/3,N]$, one
may also consider, up to rounding of the endpoints,
\[
 (N/9,2N/9]\cup(4N/9,N]
\]
and the set consisting of all odd integers together with the even
integers in $(2N/3,N]$. Each of these has asymptotic density $2/3$.
The variety of these examples suggests that the appearance of the
constant $2/3$ is not tied to a single extremal construction.

Our main result concerns cyclic groups. Although the bound
$2N/3$ is sharp for $\mathbb Z_N$ whenever $3\mid N$, the situation
changes when this obvious obstruction is absent. We show that, for
$m$ not divisible by $3$, the density of a $3$-cube-free
set is in fact bounded uniformly away from $2/3$.

\begin{thm}\label{mainthm}
There exists an absolute constant $\varepsilon>0$ such that, for
$m$ not divisible by $3$, every $3$-cube-free
set $A\subset \mathbb Z_m$ satisfies
\[
 |A|\leq \left(\frac23-\varepsilon\right)m.
\]
\end{thm}

The significance of Theorem~\ref{mainthm} is that the improvement
over $2/3$ is uniform in $m$. In particular, it rules out sequences of
$3$-cube-free sets in cyclic groups of orders not divisble by $3$ whose densities approach $2/3$ as $m\to\infty$. The value of $\varepsilon$ obtained by
our argument is explicit, although we did not attempt to optimize it.

Long and Wagner conjectured that $3$-cube-free subsets of $\mathbb Z_{2^n}$ have size at most \mbox{$\frac{5}{8} 2^n$}. We prove the following result, which in particular resolves a special case of this conjecture. 
\begin{prop}\label{prop2}
Suppose that $A\subseteq\mathbb Z_{m}$ is $3$-cube-free, where $m$ is even and $3\nmid m$.
If the number of odd residues in $A$ is at least
$m/3$, then
$|A|\leq 5m/8$.
\end{prop}

This proposition in particular proves the conjecture of Long and Wagner
when the density of $A$ among the odd residues is at least $2/3$.
This, together with an inductive argument, reduces the conjecture
to the case where this density lies in $(5/8,2/3)$.

The remainder of the paper is organized as follows.
In Section~\ref{longwagner} we prove our generalization of a special case of the Long--Wagner conjecture,
and in Section~\ref{proofofthm1} we prove Theorem~\ref{mainthm}.

\makeatletter{\let\thefootnote\relax\footnote{*Mathematical Institute of the Serbian Academy of Sciences and Arts, Kneza Mihaila 36, Belgrade 11000, Serbia. Email: zarko.randjelovic@turing.mi.sanu.ac.rs}}

\section{A result related to the Long--Wagner conjecture}\label{longwagner}

We start by defining the $2$-adic layers of $\mathbb Z_m$.
For $0\leq k\leq v_2(m)-1$, define
\[
L_k=\{x\in\mathbb Z_m\mid 2^{k}\mid x,\ 2^{k+1}\nmid x\}.
\] and define $L_{v_2(m)}=\{x\in\mathbb Z_m\mid 2^{v_2(m)}\mid x\}$.
For a $3$-cube-free subset $A$ of $\mathbb Z_m$, we use the notation
\[
\alpha_k=\frac{|L_k\cap A|}{|L_k|}
\]
for the density of $A$ in the layer $L_k$. For convenience, define $L_{v_2(m)+1}=L_{v_2(m)}$.  In order to prove Proposition \ref{prop2}, we use the following lemma, which controls the sizes of the densities.
\begin{lem}\label{lem3}
Let $m$ be an even positive integer not divisible by three and $A\subseteq\mathbb Z_m$ be $3$-cube-free.
If $\alpha_0\geq 2/3$, then $\alpha_k+\alpha_{k+1}\leq 1$
for every $1\leq k\leq v_2(m)$.
\end{lem}

\begin{proof}
First, note that $\alpha_0>2/3$, since $3\nmid |L_0|$. We claim that there is no even $x\in\mathbb Z_m$
such that $x,2x\in A$.
Indeed, suppose that such an $x$ exists.
Consider the triples $(z,z+x,z+2x)$ for all $z\in L_0$.
Since $\alpha_0>2/3$, these triples contain, on average,
more than two elements of $A$.
Hence, there is a $z\in L_0$ such that
$z,z+x,z+2x\in A$, giving
$B_3(z,x,x)\subseteq A$, a contradiction. 

For every $1\leq k\leq v_2(m)$, the map $x\mapsto 2x$
from $L_k$ to $L_{k+1}$ is surjective. It is two-to-one for $1\leq k\leq v_2(m)-2$
and bijective for $k=v_2(m)-1$. The map is also bijective for $k=v_2(m)$, since
$L_{v_2(m)}$ is an additive subgroup of $\mathbb Z_m$ of odd order. Hence,
\[
\frac{|2(A\cap L_k)|}{|L_{k+1}|}
\geq \frac{|A\cap L_k|}{|L_k|}
=\alpha_k.
\]
Since no even $x\in\mathbb Z_m$ satisfies $x,2x\in A$,
we have
\[
2(A\cap L_k)\cap(A\cap L_{k+1})=\emptyset.
\]
Therefore, $\alpha_k+\alpha_{k+1}\leq 1$. 
\end{proof}
We can also show that $2\alpha_0+\alpha_1\leq 2$.
Indeed, consider the triples $(x,2x,3x)$ for all $x\in L_0$.
Since $3\nmid m$, every element of $L_0$ appears twice
in these triples, while every element of $L_1$ appears
$|L_0|/|L_1|$ times (twice if $4\mid m$ and once otherwise).
Since $A$ is $3$-cube-free, each triple contains at most
two entries in $A$. Thus,
\[
2\alpha_0|L_0|+\alpha_1|L_0|
\leq 2|L_0|,
\]
giving $2\alpha_0+\alpha_1\leq 2$. We are now ready to prove Proposition \ref{prop2}.
\begin{proof}[Proof of Proposition \ref{prop2}]
There are fundamentally two cases, depending on how large $\alpha_0$ is. We know that $\alpha_0>2/3$.

\textit{Case 1:} $\alpha_0>3/4$. If $v_2(m)=1$, then, by the inequality
$2\alpha_0+\alpha_1\leq 2$, we have
$\alpha_1<1/2$. Hence,
\[
\frac{|A|}{m}
=\frac{\alpha_0+\alpha_1}{2}
=\frac{2\alpha_0+2\alpha_1}{4}
< \frac{2\alpha_0+\alpha_1+1/2}{4}
\leq \frac{2+1/2}{4}
=\frac58.
\]
Now suppose that $4\mid m$. If there are $x,y,x+y\in A\cap 2\mathbb Z_m$, then by averaging over all $z\in L_0$, one of the quadruples
$(z,z+x,z+y,z+x+y)$ contains all four of its elements in $A$.
Thus, $A\cap 2\mathbb Z_m$ is $2$-cube-free. Note that
\[
\bigcup_{k=2}^{v_2(m)}L_k=4\mathbb Z_m,
\]
which is a subgroup of $\mathbb Z_m$ isomorphic to
$\mathbb Z_{m/4}$. Recall that every $2$-cube-free subset of a finite abelian group $G$ has size at most $|G|/2$.
Indeed, if $B\subseteq G$ is such a set and $a\in B$, then $B$ and $a+B$
are disjoint, giving $2|B|\leq |G|$. This means that $|A\cap 4\mathbb Z_m|\leq m/8$. But we also have
\[
|A\cap L_0|+|A\cap L_1|
=\alpha_0|L_0|+\alpha_1|L_1|
=(2\alpha_0+\alpha_1)|L_1|
\leq 2|L_1|=\frac{m}{2}.
\]
Together with $|A\cap 4\mathbb Z_m|\leq m/8$, this gives
\[
|A|\leq \frac{m}{2}+\frac{m}{8}=\frac{5m}{8},
\]
which completes Case 1.

\textit{Case 2:} $2/3<\alpha_0\leq 3/4$. If $v_2(m)=1$ then by Lemma \ref{lem3} $\alpha_1\le 1/2$. The rest is done in the same way as in Case 1. Thus, we may assume that $v_2(m)>1$. Let $r=v_2(m)$. Note that
\[
\frac{|A|}{m}
=\frac{\alpha_0}{2}+\frac{\alpha_1}{4}
+\cdots+\frac{\alpha_{r-1}}{2^r}
+\frac{\alpha_r}{2^r}.
\]Here we will use all the inequalities on the $\alpha_k$ and show that they are enough. In view of this, we prove the following claim.

\medskip
\noindent\textit{Claim.}
For $\ell\ge 2$, let $\beta_0,\beta_1,\ldots,\beta_{\ell}$ be real numbers
such that $\beta_\ell\le 1/2$, $\beta_0\le 3/4$, 
$2\beta_0+\beta_1\leq 2$, and
$\beta_k+\beta_{k+1}\leq 1$ for every $1\leq k\leq \ell-1$.
Then
\[
\sum_{k=0}^{\ell-1}\frac{\beta_k}{2^k}
+\frac{\beta_\ell}{2^{\ell-1}}\leq\frac{5}{4}.
\]
\medskip
\noindent\textit{Proof of the claim.} We will prove the inequality by expressing its left-hand side as a nonnegative linear combination of the left-hand sides of the assumed inequalities. Indeed,
\[
\sum_{k=0}^{\ell}\frac{\beta_k}{2^k}
+\frac{\beta_\ell}{2^\ell}
=\frac13\beta_0+\frac13(2\beta_0+\beta_1)
+\sum_{k=1}^{\ell-1}\frac{\beta_k+\beta_{k+1}}{3\cdot 2^k}
+\frac{\beta_\ell}{3\cdot 2^{\ell-2}}.
\] Thus, the left-hand side is at most
\[
\begin{aligned}
\frac{1}{3}\cdot\frac{3}{4}+\frac{1}{3}\cdot 2
+\frac{1}{3}\left(\sum_{k=1}^{\ell-1}\frac{1}{2^k}
+\frac{1}{2^{\ell-1}}\right)
=\frac{1}{4}+\frac{2}{3}+\frac{1}{3}
=\frac{5}{4}.
\end{aligned}
\]

By Lemma~\ref{lem3}, we have
$\alpha_{v_2(m)}\leq 1/2$, since
$L_{v_2(m)+1}=L_{v_2(m)}$. Applying the claim with $\beta_k=\alpha_k$ and
$\ell=v_2(m)$, we obtain the desired bound.
\end{proof}

\section{Proof of Theorem 1}\label{proofofthm1}
We begin by outlining the idea of the proof.
Let $A\subseteq\mathbb Z_m$ be a $3$-cube-free set, and for any $x\in\mathbb Z_m$ define
$T_x=(x,2x,3x)$. Note that no $T_x$ can have all its entries in $A$. We say that $T_x$ is \emph{good} if two of its entries are in $A$. The idea is to construct a fixed set $B$ (independent of $m$) such that there is no $x$ for which $T_{bx}$ is good for every $b\in B$. Afterward, we will show that if the density of $A$ is sufficiently close to $2/3$ there is indeed such an $x$, finishing the proof of the theorem. In view of this, we prove the following lemma.
\begin{lem}\label{Blemma}
Let $C=\{1,2,3,4,5,7,9\}$ and $B=C\cup 3C\cup 9C$.
Then there is no $x\in\mathbb Z_m$ such that $T_{bx}$ is good
for every $b\in B$.
\end{lem}
\begin{proof}
Suppose that there is an $x\in\mathbb Z_m$ such that $T_{bx}$ is good for every $b\in B$. We first claim that there is a $y\in\{x,3x,9x\}$ such that $y,3y\in A$. 

Supposing the contrary and using that $T_x$, $T_{3x}$, and $T_{9x}$ are good,
we have that in the sequence $x,3x,9x,27x$, no two consecutive elements
are both in or both outside of $A$. Thus, there is a $z\in\{x,3x\}$ such that $z,9z\notin A$. Since $T_z$ is good, $2z,3z\in A$, and since $T_{3z}$ and $T_{9z}$ are good, $6z,18z\in A$ as well. Now $4z,12z\notin A$, contradicting that $T_{4z}$ is good.

Note that $T_{cy}$ is good for every $c\in C$. Since $T_y$ is good, $2y\notin A$, and since $T_{2y}$ is good, $4y,6y\in A$. Furthermore, since $T_{3y}$ is good, $9y\notin A$, and since $T_{9y}$ is good, $18y\in A$. Now we get that $12y\notin A$, and since $T_{4y}$ is good, $8y\in A$. Since $B_3(y,3y,3y)\nsubseteq A$, we have $7y\notin A$, and since $T_{7y}$ is good, $14y,21y\in A$. By considering $B_3(4y,4y,6y)$, we obtain that $10y\notin A$,
and since $T_{5y}$ is good, $15y\in A$. Now we have that $B_3(3y,3y,15y)\subseteq A$, giving a contradiction.

\end{proof}

We are now ready to prove Theorem~\ref{mainthm}.

\begin{proof}[Proof of Theorem~\ref{mainthm}] Assume that $\varepsilon$ is sufficiently small (we will choose $\varepsilon$ later) and let $A$ be a $3$-cube-free subset of $\mathbb Z_m$ of size greater than $\left(\frac{2}{3}-\varepsilon\right)m$. For convenience, if $v_2(m)<3$, let $L_k=L_{v_2(m)}$ for $v_2(m)< k\le 3$. We will ensure that $\varepsilon<1/192$, which guarantees
\[
\frac{2}{3}-\varepsilon
>\frac{2}{3}\left(1-\frac{1}{8}\right)
+\frac{5}{8}\cdot\frac{1}{8}.
\] 
We can inductively show that $\alpha_k\leq 2/3$ for $k\leq 3$. Indeed, if $k\geq v_2(m)$, then
$L_k=2^{v_2(m)}\mathbb Z_m$ is a subgroup of
$\mathbb Z_m$, so $\alpha_k\leq 2/3$ by Meng's bound. If $k<v_2(m)$ and $\alpha_k>2/3$, then by Proposition \ref{prop2},
applied to $2^k\mathbb Z_m$, we have
\[
|A\cap 2^k\mathbb Z_m|\leq \frac{5}{8}\cdot\frac{m}{2^k}.
\]
Hence, by the induction hypothesis,
\[
\frac{|A|}{m}
\leq \frac{2}{3}\left(1-\frac{1}{2^k}\right)
+\frac{5}{8}\cdot\frac{1}{2^k}
\leq \frac{2}{3}-\frac{1}{192}
<\frac{2}{3}-\varepsilon.
\]
We will show that for any $\delta>0$, if $\varepsilon$ is sufficiently small,
then $\alpha_k\ge 2/3-\delta$ for every $0\leq k\leq 3$.
Indeed, it is enough to consider
$0\leq k\leq \min\{3,v_2(m)\}$.
Suppose that $\alpha_k< 2/3-\delta$.
We already have $\alpha_i\leq 2/3$ for every $i<k$.
Moreover, if $k<v_2(m)$, the remaining subgroup
$2^{k+1}\mathbb Z_m$ has density at most $2/3$ by Meng's bound.
When $k=v_2(m)$, there is no remaining subgroup, and the same argument works.
Consequently,
\[
\frac{|A|}{m}
< \frac{2}{3}-\delta\frac{|L_k|}{m}
\leq \frac{2}{3}-\frac{\delta}{16},
\]
giving a contradiction if $\varepsilon<\delta/16$. Now assume that $\alpha_k\geq 2/3-\delta$ for every $0\leq k\leq 3$,
where $\delta>0$ is sufficiently small.
For each $0\leq k\leq 2$, consider the triples $T_x$ for $x\in L_k$.
Since $3\nmid m$, both $x$ and $3x$ vary uniformly over $L_k$,
while $2x$ varies uniformly over $L_{k+1}$.
Hence, the total number of appearances of elements of $A$ in these triples
is at least
\[
(2\alpha_k+\alpha_{k+1})|L_k|
\geq (2-3\delta)|L_k|.
\]
On the other hand, if $x$ is chosen uniformly at random from $L_k$,
the total number of entries of $T_x$ belonging to $A$ is at most
\[
2\mathbb P(T_x\text{ is good})|L_k|
+\mathbb P(T_x\text{ is not good})|L_k|
=2|L_k|-\mathbb P(T_x\text{ is not good})|L_k|.
\]
Therefore,
\[
\mathbb P(T_x\text{ is not good})\leq 3\delta.
\]
Let $B$ be as in Lemma~\ref{Blemma} and let
\begin{align}\label{eq1}
\delta<\frac{1}{3\sum_{b\in B}b^*},
\end{align}
where $b^*$ is the largest divisor of $b$ coprime to $6$. For any $b\in B$, note that as $x$ varies uniformly over $L_0$,
$bx$ varies uniformly over $b^*L_{v_2(b)}$, since multiplication
by $3$ only permutes the layers.
We also have that
\[
|b^*L_{v_2(b)}|\geq \frac{|L_{v_2(b)}|}{b^*}.
\]
This means that the probability that $T_{bx}$ is not good,
when $x$ is chosen uniformly at random from $L_0$,
is at most $3\delta b^*$. Thus, by \eqref{eq1}, there is an $x\in L_0$ such that $T_{bx}$ is good for every $b\in B$, contradicting Lemma~\ref{Blemma}. This completes the proof. Since $\sum_{b\in B}b^*=47$, we can take
$\varepsilon=1/2257$ and choose
\[
\delta\in\left(\frac{16}{2257},\frac{1}{141}\right).
\]
\end{proof}

\bibliographystyle{plain}
\bibliography{3cube}

@article{LongWagner2019,
author = {Long, Jason and Wagner, Adam Zsolt},
title = {The largest projective cube-free subsets of {$\mathbb{Z}_{2^n}$}},
year = {2019},
issue_date = {Oct 2019},
publisher = {Academic Press Ltd.},
address = {GBR},
volume = {81},
number = {C},
issn = {0195-6698},
url = {https://doi.org/10.1016/j.ejc.2019.05.005},
doi = {10.1016/j.ejc.2019.05.005},
journal = {Eur. J. Comb.},
month = oct,
pages = {156–171},
numpages = {16}
}

@article{Rado1933,
  author  = {Richard Rado},
  title   = {Studien zur Kombinatorik},
  journal = {Math. Z.},
  volume  = {36},
  year    = {1933},
  pages   = {424--480}
}

@article{Meng2026, 
title={On Cube-Free Problems}, 
volume={33}, 
url={https://www.combinatorics.org/ojs/index.php/eljc/article/view/v33i1p16}, DOI={10.37236/14052}, 
number={1}, 
journal={The Electronic Journal of Combinatorics}, 
author={Meng, Yuchen},
year={2026}, 
month={Jan.}, 
pages={P1.16} 
}

@article{GreenRuzsa2005,
  author  = {Ben Green and Imre Z. Ruzsa},
  title   = {Sum-free sets in abelian groups},
  journal = {Israel Journal of Mathematics},
  volume  = {147},
  pages   = {157--188},
  year    = {2005},
  doi     = {10.1007/BF02785363}
}

@article{EberhardMellon2025,
  author  = {Sean Eberhard and Connor Mellon},
  title   = {Normal covering numbers for {$S_n$} and {$A_n$}
             and additive combinatorics},
  journal = {Bulletin of the London Mathematical Society},
  volume  = {57},
  number  = {11},
  pages   = {3307--3325},
  year    = {2025},
  doi     = {10.1112/blms.70154}
}

\end{document}